# A class of unbiased location invariant Hill-type estimators for heavy tailed distributions

**Jiaona Li[1], Zuoxiang Peng[1] and Saralees Nadarajah[2]**

[1]*School of Mathematics and Statistics, Southwest University, Chongqing, 400715*

[2]*School of Mathematics, University of Manchester, Manchester, United Kingdom*

**Abstract:** Based on the methods provided in Caeiro and Gomes (2002) and Fraga Alves (2001), a new class of location invariant Hill-type estimators is derived in this paper. Its asymptotic distributional representation and asymptotic normality are presented, and the optimal choice of sample fraction by Mean Squared Error is also discussed for some special cases. Finally comparison studies are provided for some familiar models by Monte Carlo simulations.

## 1. Introduction

In this paper we study the asymptotic properties of a class of unbiased location invariant Hill-type heavy tailed index estimators. Let $\{X_n, n \geq 1\}$ be an independent and identically distributed (i.i.d) sequence with distribution function (d.f) $F(x)$ belonging to the domain of attraction of an extreme value distribution,

$$G_\gamma(x) = \exp\left\{-(1+\gamma x)^{-1/\gamma}\right\}, \ \gamma \in \mathbb{R} \ \text{and} \ 1+\gamma x > 0,$$

i.e. there exist normalizing constants $a_n > 0$ and $b_n \in \mathbb{R}$ such that

$$F^n(a_n x + b_n) \longrightarrow G_\gamma(x) \tag{1.1}$$

as $n \to \infty$ for all $x \in \mathbb{R}$. We write $F \in D(G_\gamma)$ if (1.1) holds and call $\gamma$ the tail index.

The applications of the heavy tailed distributions may be found in many fields such as insurance, finance, climatology and environmental science (cf. Embrechets et al. (1997)), and the problem of estimating the heavy tailed index has been studied extensively. For positive $\gamma$, Hill (1975) introduced the well known Hill estimator, which is

$$\widehat{\gamma}_{n,H}(k) = \frac{1}{k}\sum_{i=1}^{k-1} \ln X_{n-i,n} - \ln X_{n-k,n},$$

where $X_{1,n} \leq X_{2,n} \leq \cdots \leq X_{n,n}$ are the order statistics of $X_1, X_2, \cdots, X_n$. As an extension of Hill-type estimator for general $\gamma \in \mathbb{R}$, Dekkers et al. (1989) proposed the following moment estimator

$$\widehat{\gamma}_n^M(k) = M_n^{(1)} + 1 - \frac{1}{2}\left(1 - \frac{(M_n^{(1)})^2}{M_n^{(2)}}\right)^{-1},$$

where

$$M_n^{(j)}(k) = \frac{1}{k}\sum_{i=0}^{k-1}(\ln X_{n-i,n} - \ln X_{n-k,n})^j, \ j = 1, 2.$$

In applications both Hill and moment estimators are sensitive to the threshold, say $X_{n-k,n}$, and these two estimators are not invariant to the affine transformation. Fraga Alves (2001) established the scale and location invariant Hill-type estimator given by

$$\widehat{\gamma}_n^H(k_0, k) = \frac{1}{k_0}\sum_{i=0}^{k_0-1}\ln\frac{X_{n-i,n} - X_{n-k,n}}{X_{n-k_0,n} - X_{n-k,n}},$$

where the intermediate sequences $k_0$ and $k$ satisfy

$$k = k_n = o(n), \ k_0 = o(k_n), \ k_n \to \infty, \ k_0 \to \infty, \ \text{as } n \to \infty. \tag{1.2}$$

The asymptotic properties of $\widehat{\gamma}_n(k_0, k)$ have been considered when

$$\frac{U(tx)/U(t) - x^\gamma}{a(t)} \to x^\gamma \frac{x^\rho - 1}{\rho}$$

holds as $t \to \infty$, where $U(t) = F^\leftarrow(1 - 1/t)$, $t \geq 1$, where $F^\leftarrow(x)$ denotes the inverse function of $F(x)$. Based on the methods provided in Fraga Alves (2001) and Dekkers et al. (1989), Ling et al. (2007a, 2007b) proposed a kind of location invariant moment-type tail index estimator and considered its asymptotic properties under some second order regular varying conditions.

Meanwhile the problem of searching for the unbiased estimators has been considered by many authors, see, e.g., Peng (1998), Beirlant et al. (2002), Caeiro and Gomes (2002), Gomes and Martins (2002, 2004), Peng and Qi (2006a, 2006b), Gomes et al. (2008), Gomes and Henriques (2008), and Qi (2008a, 2008b). Here we are interested in the class of semi-parametric heavy tail estimator introduced by Caeiro and Gomes (2002):

$$\gamma_n^{(\alpha)}(k) = \frac{\Gamma(\alpha)}{M_n^{(\alpha-1)}(k)}\left(\frac{M_n^{(2\alpha)}(k)}{\Gamma(2\alpha+1)}\right)^{1/2}, \ \alpha \geq 1,$$

where $\Gamma(\cdot)$ is the gamma function and

$$M_n^{(\alpha)}(k) := \frac{1}{k}\sum_{i=1}^{k-1}(\ln X_{n-i,n} - \ln X_{n-k,n})^\alpha, \ \alpha > 0.$$

They considered the asymptotic distributional representation of $M_n^{(\alpha)}(k)$ and the choice of tuning parameter $\alpha$ such that $\gamma_n^{(\alpha)}(k)$ is asymptotically normal with asymptotic null bias under the assumptions that

$$\frac{\ln U(tx) - \ln U(t) - \gamma \ln x}{A(t)} \to \frac{x^\rho - 1}{\rho}, \quad x > 0$$

and $\sqrt{k}A(n/k) \to \lambda$, finite.

Motivated by the works of Fraga Alves (2001) and Caeiro and Gomes (2002), we propose a new class of location invariant estimators for a heavy tailed distribution based on the asymptotic distributional representation of the following statistic:

$$M_n^{(\alpha)}(k_0, k) = \frac{1}{k_0} \sum_{i=0}^{k_0-1} \left\{ \ln \frac{X_{n-i,n} - X_{n-k,n}}{X_{n-k_0,n} - X_{n-k,n}} \right\}^\alpha, \quad \alpha \in \mathbb{R}^+. \tag{1.3}$$

The asymptotic distributional representation of $M_n^{(\alpha)}(k_0, k)$ will be derived under the following second order regular variation condition:

$$\lim_{t \to \infty} \frac{\frac{U(tx) - U(t)}{a(t)} - \frac{x^\gamma - 1}{\gamma}}{A(t)} = \Psi_{\gamma,\rho}(x) \quad \text{for all } x > 0, \tag{1.4}$$

where

$$\Psi_{\gamma,\rho}(x) = \begin{cases} \dfrac{x^{\gamma+\rho} - 1}{\gamma + \rho}, & \text{if } \gamma + \rho \neq 0, \\ \ln x, & \text{if } \gamma + \rho = 0. \end{cases}$$

and $\rho < 0$ is the second order parameter and $|A(t)| \in RV_\rho$ (cf. Corollary 2.3.5 of de Haan and Ferreira (2006)). Here, $f \in RV_\beta$ means $\lim_{t \to \infty} f(tx)/f(t) = x^\beta$ for all $x > 0$.

Based on the convergence $M_n^{(\alpha)}(k_0, k) \xrightarrow{p} \Gamma(\alpha+1)\gamma^\alpha$ as $k_0 \to \infty$, $k_0 = o(k)$, a location invariant Hill-type estimator for the heavy tailed index may be defined by

$$\widehat{\gamma}_n^{(\alpha)}(k_0, k) = \frac{\Gamma(\alpha)}{M_n^{(\alpha-1)}(k_0, k)} \left[ \frac{M_n^{(2\alpha)}(k_0, k)}{\Gamma(2\alpha+1)} \right]^{\frac{1}{2}}, \quad \alpha \geq 1, \tag{1.5}$$

which converges to $\gamma$ in probability. Finite sample simulation shows that the positive tuning parameter $\alpha$ may be chosen appropriately to improve the performance of tail index estimation in applications.

The paper is organized as follows. In section 2, we provide the main results, i.e. the asymptotic distributional representations of (1.3) and (1.5), and the optimal choice of the sample fraction $k_0$ by mean squared error (MSE) for some special distributions. Related proofs are deferred to section 4. Simulation studies are performed in section 3.

## 2. The main results

Throughout this paper, we assume that $X_1, X_2, \cdots, X_n$ are i.i.d random variables (r.v.s) with d.f $F(x)$, and denote by $X_{1,n} \leq X_{2,n} \leq \cdots \leq X_{n,n}$ the order statistics of $X_1, \cdots, X_n$. For the heavy tail distribution, i.e. $\gamma > 0$, we know that

$$F \in D(G_\gamma) \Leftrightarrow 1 - F \in RV_{-1/\gamma} \Leftrightarrow U \in RV_\gamma. \tag{2.1}$$

We need the following notations to simplify the statements of our main results:

$$\sigma_\alpha = \sqrt{\Gamma(2\alpha+1) - \Gamma^2(\alpha+1)},$$

$$\mu_\alpha(\rho) = \frac{\Gamma(\alpha)}{\rho} \frac{1 - (1-\rho)^\alpha}{(1-\rho)^\alpha},$$

$$V_\alpha = \frac{1}{4} \left\{ \frac{\Gamma(4\alpha)}{\alpha \Gamma^2(2\alpha)} + \frac{4\Gamma(2\alpha-1)}{\Gamma^2(\alpha)} - \frac{2\Gamma(3\alpha)}{\alpha \Gamma(\alpha)\Gamma(2\alpha)} - 1 \right\} \tag{2.2}$$

and

$$b_\alpha(\gamma) = (1+\gamma)^{1-\alpha} - \frac{1}{2}(1+\gamma)^{-2\alpha} - \frac{1}{2}. \tag{2.3}$$

The first result is about the asymptotic distributional representation of $M_n^{(\alpha)}(k_0, k)$.

**Theorem 2.1.** *Suppose that (2.1) holds for $\gamma > 0$, and the intermediate $k_0$ and $k$ satisfy (1.2). Then $M_n^{(\alpha)}(k_0, k)$ converges in probability to $\Gamma(\alpha+1)\gamma^\alpha$. Furthermore, if the second order framework in (1.4) holds, we may obtain the following asymptotic distributional representation*

$$M_n^{(\alpha)}(k_0, k) \stackrel{d}{=} \gamma^\alpha \Gamma(\alpha+1) + \frac{\gamma^\alpha \sigma_\alpha}{\sqrt{k_0}} P_n^{(\alpha)} + \alpha \gamma^\alpha \mu_\alpha(-\gamma) \left(\frac{k_0}{k}\right)^\gamma (1 + o_P(1)) + B_n,$$

*where $P_n^{(\alpha)}$ is an asymptotically standard normal r.v., and*

$$B_n = \begin{cases} \dfrac{\alpha \gamma^\alpha \rho \mu_\alpha(\rho)}{\gamma + \rho} A(n/k) \left(\dfrac{k_0}{k}\right)^{-\rho} (1 + o_P(1)), & \text{if } \gamma + \rho \neq 0, \\ \alpha \gamma^{\alpha+1} \mu_\alpha(-\gamma) A(n/k) \left(\dfrac{k_0}{k}\right)^\gamma \ln\left(\dfrac{k_0}{k}\right)(1 + o_P(1)), & \text{if } \gamma + \rho = 0. \end{cases}$$

Based on Theorem 2.1, we may derive the asymptotic distributional representation of the proposed estimator $\widehat{\gamma}_n^{(\alpha)}(k_0, k)$ in (1.5), which is the following result.

*J. Li et al./Unbiased location invariant Hill-type estimators* 833**Theorem 2.2.** *Under the conditions of Theorem 2.1, the following asymptotic distributional representation*

$$\widehat{\gamma}_n^{(\alpha)}(k_0, k) \stackrel{d}{=} \gamma + \frac{\gamma\sqrt{V_\alpha}}{\sqrt{k_0}} T_n^{(\alpha)} + b_\alpha(\gamma)\left(\frac{k_0}{k}\right)^\gamma + o_P\left(\frac{1}{\sqrt{k_0}}\right) + o_P\left(\left(\frac{k_0}{k}\right)^\gamma\right) + R_n \qquad (2.4)$$

*holds, where $T_n^{(\alpha)}$ is asymptotically standard normal, and*

$$R_n = \begin{cases} -\dfrac{\gamma}{\gamma+\rho} b_\alpha(-\rho) A(n/k)\left(\dfrac{k_0}{k}\right)^{-\rho}(1+o_P(1)), & \text{if } \gamma+\rho \neq 0, \\ \gamma b_\alpha(\gamma) A(n/k)\left(\dfrac{k_0}{k}\right)^\gamma \ln\left(\dfrac{k_0}{k}\right)(1+o_P(1)), & \text{if } \gamma+\rho = 0. \end{cases}$$

*Furthermore, if there exist finite $\lambda_1$ and $\lambda_2$ such that*

$$\sqrt{k_0}\left(\frac{k_0}{k}\right)^\gamma \to \lambda_1 \text{ and } \sqrt{k_0} A(n/k) \to \lambda_2,$$

*we have*

$$\sqrt{k_0}\left(\widehat{\gamma}_n^{(\alpha)}(k_0,k) - \gamma\right) \to N\left(\lambda_1 b_\alpha(\gamma), \gamma^2 V_\alpha\right),$$

*where $V_\alpha$ and $b_\alpha(\gamma)$ are defined in (2.2) and (2.3), respectively. Consequently, for every $\gamma > 0$, there exists $\alpha_0$ given by*

$$\alpha_0 = \alpha_0(\gamma) = \frac{\ln\left(1+\gamma+\sqrt{(1+\gamma)^2 - 1}\right)}{\ln(1+\gamma)}$$

*such that $b_{\alpha_0}(\gamma) = 0$, i.e. $\widehat{\gamma}_n^{(\alpha_0)}(k_0, k)$ has asymptotic null bias, even when $\sqrt{k_0}(k_0/k)^\gamma \to \lambda_1 \neq 0$ and $\sqrt{k_0} A(n/k) \to \lambda_2 \neq 0$.*

For special $A(t)$, we will consider the optimal choice of the sample fraction $k_0$ as a function of $k$, $\gamma$, $\rho$ and $\alpha$ following the criterion of Fraga Alves (2001) to compute $k_0 \equiv k_0^{opt}$ such that the asymptotic MSE (Fraga Alves, 2001) of $\widehat{\gamma}_n^{(\alpha)}(k_0, k)$ denoted by $MSE_\infty(\widehat{\gamma}_n^{(\alpha)}(k_0, k))$ is minimal. Note that an asymptotic MSE is just the usual MSE based on some asymptotic relationship. For example, suppose $\widehat{\theta}_n$ is an estimator of $\theta$ based on a random sample of size $n$ and that it satisfies $\sqrt{a_n}(\widehat{\theta}_n - \theta - b_n) \to N(\mu, \sigma^2)$ as $n \to \infty$ for some $a_n$, $b_n$, $\mu$ and $\sigma^2$. Then the asymptotic MSE of $\widehat{\theta}_n$ denoted by $MSE_\infty(\widehat{\theta}_n)$ is simply $\{b_n + \mu/\sqrt{a_n}\}^2 + \sigma^2/a_n$.

**Theorem 2.3.** *Suppose that (1.4) holds with $A(t) \sim ct^\rho$, $\rho < 0$ and $c \neq 0$. Set*

$$k_0^{(1)} = \left[\frac{\gamma V_\alpha}{2 b_\alpha^2(\gamma)}\right]^{\frac{1}{2\gamma+1}} \cdot k^{\frac{2\gamma}{2\gamma+1}};$$

$$k_0^{(2)} = \left[\frac{(\gamma+\rho)^2 V_\alpha}{-2c^2 \rho b_\alpha^2(-\rho)}\right]^{\frac{1}{-2\rho+1}} \cdot n^{\frac{-2\rho}{-2\rho+1}};$$

$$k_0^{(3)} = \left[\frac{c\gamma b_\alpha(-\rho)}{(\gamma+\rho) b_\alpha(\gamma)}\right]^{\frac{1}{\gamma+\rho}} \cdot k^{\frac{\gamma}{\gamma+\rho}} \cdot n^{\frac{\rho}{\gamma+\rho}},$$

where $V_\alpha$ and $b_\alpha(\gamma)$ are defined as before. Then we may obtain the sequence $k_0^{opt}$ that minimizes $MSE_\infty(\widehat{\gamma}_n^{(\alpha)}(k_0, k))$,

  i. For $\gamma \leq -\rho$, $k_0^{opt} = k_0^{(1)}$;
  
  ii. For $\gamma > -\rho$,
  
   (a). If $k \ll n^{-\rho(2\gamma+1)/\{\gamma(-2\rho+1)\}}$, $k_0^{opt} = k_0^{(1)}$;
   
   (b). If $k \gg n^{-\rho(2\gamma+1)/\{\gamma(-2\rho+1)\}}$, $k_0^{opt} = k_0^{(2)}$ if $cb_\alpha(-\rho)b_\alpha(\gamma) < 0$ and $k_0^{opt} = k_0^{(3)}$ if $cb_\alpha(-\rho)b_\alpha(\gamma) > 0$;
   
   (c). If $k \sim Dn^{-\rho(2\gamma+1)/\{\gamma(-2\rho+1)\}}$ with $D \neq 0$, then $k_0^{opt} \sim D_1 n^{-2\rho/(-2\rho+1)}$ with $D_1 \equiv D_1(\gamma, \rho, \alpha, c, D)$ satisfying

$$a_1 D_1^{2\gamma+1} + a_2 D_1^{\gamma-\rho+1} + a_3 D_1^{-2\rho+1} = \gamma^2 V_\alpha,$$

where $a_1 = 2\gamma b_\alpha^2(\gamma)D^{-2\gamma}$, $a_2 = [2c\gamma(\rho-\gamma)/(\gamma+\rho)]b_\alpha(\gamma)b_\alpha(-\rho)D^{-\gamma}$ and $a_3 = -2\rho[c\gamma b_\alpha(-\rho)/(\gamma+\rho)]^2$.

**Corollary 2.1.** *Suppose (1.4) holds for $A(t) \sim ct^\rho$ with $\rho < 0$, $c \neq 0$, and the intermediate sequence $k$ satisfies $k \to \infty$, $k/n \to 0$ as $n \to \infty$ for $\gamma \leqslant -\rho$; or $k \to \infty$, $kn^{-\rho(2\gamma+1)/\{\gamma(-2\rho+1)\}} \to 0$ for $\gamma > -\rho$. Choosing $k_0 \sim k_0^{(1)} = [\gamma V_\alpha b_\alpha^{-2}(\gamma)/2]^{1/(2\gamma+1)} \cdot k^{2\gamma/(2\gamma+1)}$, we have*

$$\sqrt{k_0^{(1)}} \left( \widehat{\gamma}_n^{(\alpha)}(k_0^{(1)}, k) - \gamma \right) \to N\left(\lambda_1 b_\alpha(\gamma), \gamma^2 V_\alpha\right)$$

*with $b_\alpha(\gamma)$ and $V_\alpha$ defined as before.*

## 3. Simulation study

Firstly, we present some $\alpha_0(\gamma)$ for heavy index $\gamma$ such that $b_{\alpha_0}(\gamma) = 0$.

Table 1 shows that $\alpha_0(\gamma)$ is a decreasing function of $\gamma$. Firstly we consider the effect of the tuning parameter $\alpha$ on the heavy tail index estimator proposed in this paper. We randomly select a sample from Fréchet d.f $F(x) = \exp(-x^{-1/\gamma})$ with $\gamma = 1$. For $\gamma = 1$, choose $\alpha_0 = 1.90$ such that $b_{\alpha_0}(\gamma) = 0$ (cf. Table 1). We compare the sample paths of $\widehat{\gamma}_n^{(1)}$, $\widehat{\gamma}_n^{(1.90)}$ and $\widehat{\gamma}_n^{(3)}$ with sample size $n = 3000$. Figure 1 shows that $\widehat{\gamma}_n^{(1.90)}$ has a much smaller bias than others.

For the optimal sequence $k_0$, we denote $\widetilde{k}_0 := \arg\min_{k_0} MSE(\widetilde{\gamma}_n^H(k_0, k))$ and $\widehat{k}_0 := \arg\min_{k_0} MSE(\widehat{\gamma}_n^{(\alpha)}(k_0, k))$. Consequently, we get the following possible measure of asymptotic relative efficiency (AREFF) of the proposed location invariant estimator $\widehat{\gamma}_n^{(\alpha)} = \widehat{\gamma}_n^{(\alpha)}(\widehat{k}_0(\alpha), k)$ to the Hill-type location invariant estimator $\widetilde{\gamma}_n^H = \widetilde{\gamma}_n^H(k_0, k)$ (bias-corrected).

$$AREFF_{[\widehat{\gamma}_n^{(\alpha)}, \widetilde{\gamma}_n^H]} = V_\alpha^{-\gamma/(2\gamma+1)} \cdot [\gamma/\{(1+\gamma)|b_\alpha(\gamma)|\}]^{1/(2\gamma+1)} \quad (3.1)$$

whenever $b_\alpha(\gamma) \neq 0$.

Table 1
$\alpha_0(\gamma)$ as a function of $\gamma$, i.e. $\alpha_0(\gamma) = \{\alpha : b_\alpha(\gamma) = 0\}$.

| $\gamma$ | 0 | 0.1 | 0.5 | 0.75 | 1 | 1.5 | 2 | 2.5 | 3 | 4 | $\infty$ |
|---|---|---|---|---|---|---|---|---|---|---|---|
| $\alpha_0(\gamma)$ | $\infty$ | 4.65 | 2.37 | 2.07 | 1.9 | 1.71 | 1.60 | 1.54 | 1.49 | 1.42 | 1 |

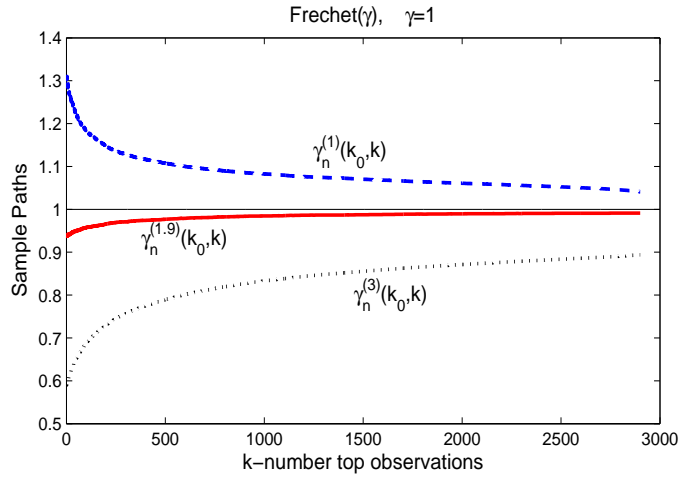

FIG 1. *Sample paths of $\widehat{\gamma}_n^{(1)}(k_0, k)$, $\widehat{\gamma}_n^{(1.90)}(k_0, k)$ and $\widehat{\gamma}_n^{(3)}(k_0, k)$ for Fréchet (1) distribution with sample size $n = 3000$.*

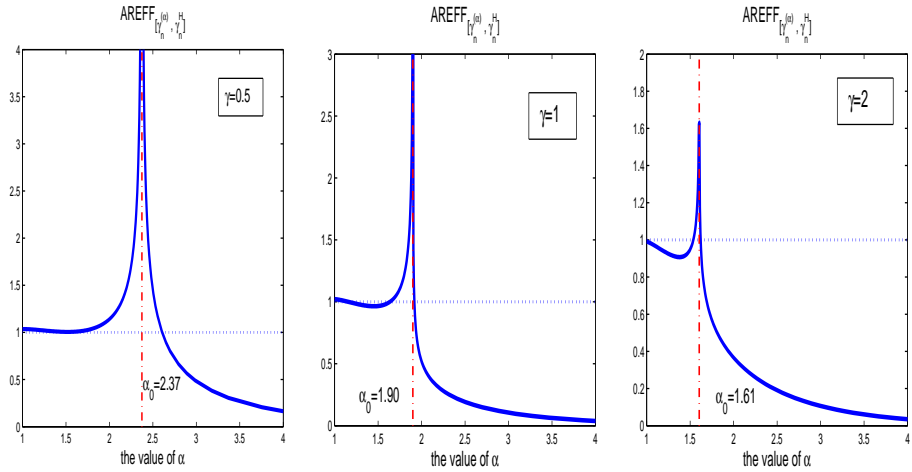

FIG 2. *Asymptotic efficiency of $\widehat{\gamma}_n^{(\alpha)}$ relative to $\widetilde{\gamma}_n^H$ for $\gamma = 0.5$, 1 and 2.*

Figure 2 shows the $AREFF_{[\widehat{\gamma}_n^{(\alpha)}, \widetilde{\gamma}_n^H]}$ for $\gamma = 0.5$, 1 and 2. Simulation shows that for every $\gamma$, we may find some $\alpha$ on the left region of $\alpha_0$ such that the AREFF in (3.1) is greater than one.

Next we compare the relative efficiency of the proposed location invariant Hill-type estimator and that of Fraga Alves (2001) in terms of average mean and MSE for finite sample size. We consider the following two models:

– Burr $(\alpha, \beta)$ distribution with d.f $F(x) = 1 - (1 + x^\alpha)^{-\beta}$, where $x \geq 0$, $\alpha > 0$ and $\beta > 0$.

– Pareto $(\gamma)$ distribution with d.f $F(x) = 1 - x^{-1/\gamma}$, where $x > 0$ and $\gamma > 0$.

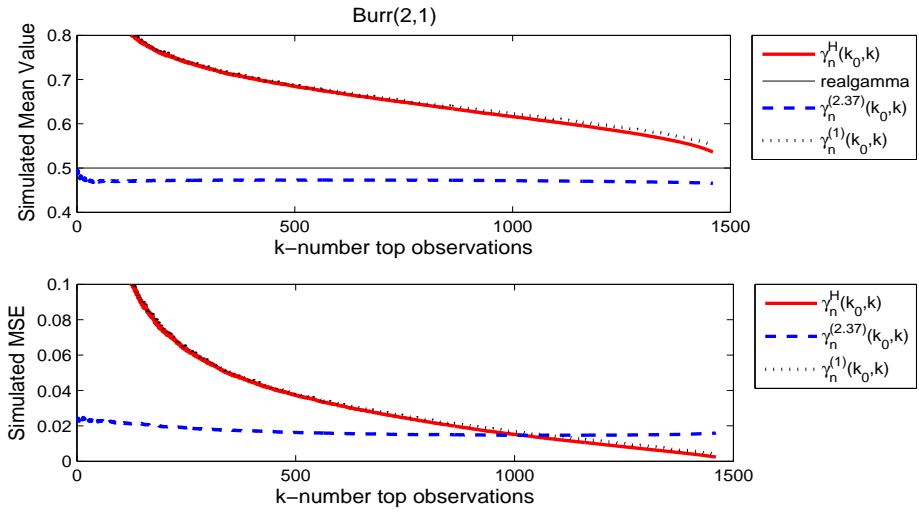

Fig 3. *Mean value and MSE of $\widetilde{\gamma}_n^H(k_0,k)$, $\widehat{\gamma}_n^{(\alpha_0)}(k_0,k)$ and $\widehat{\gamma}_n^{(1)}(k_0,k)$ for Burr (2,1) model with $\gamma = 0.5$ and sample size $n = 1500$.*

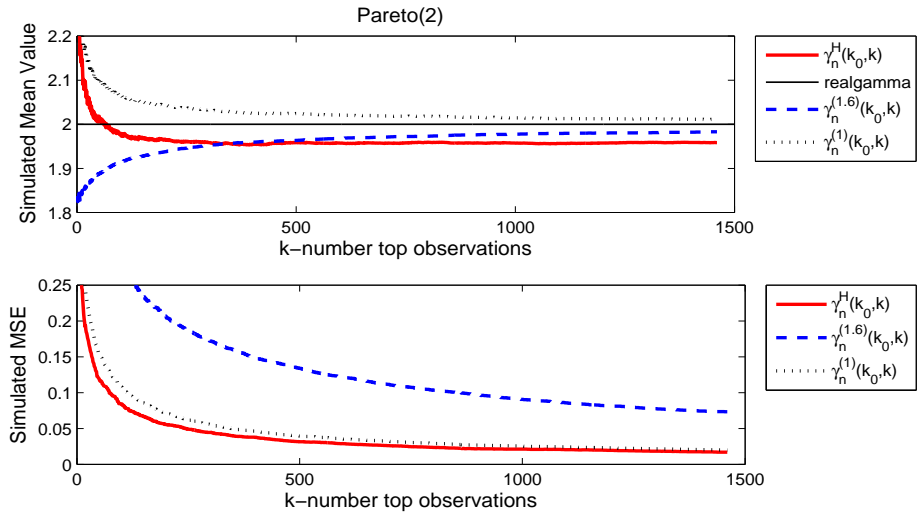

Fig 4. *Mean value and MSE of $\widetilde{\gamma}_n^H(k_0,k)$, $\widehat{\gamma}_n^{(\alpha_0)}(k_0,k)$ and $\widehat{\gamma}_n^{(1)}(k_0,k)$ for Pareto model with $\gamma = 2$ and sample size $n = 1500$.*

Figures 3 and 4 show the simulated mean value and MSE of $\widetilde{\gamma}_n^H(k_0,k)$, $\widehat{\gamma}_n^{(1)}(k_0,k)$ and $\widehat{\gamma}_n^{(\alpha_0)}(k_0,k)$ for the distributions just mentioned. For both distributions, $\widehat{\gamma}_n^{(\alpha_0)}(k_0,k)$ has a small bias and is closer to the real extreme value index value (see Figures 3 and 4).

TABLE 2
*Coverage probability comparison for 95% asymptotic confidence intervals based on* (3.2) *and* (3.3).

|  | $\widetilde{\gamma}_n^H(k_0,k)$ | $\widehat{\gamma}_n^{(2.37)}(k_0,k)$ | $\widetilde{\gamma}_n^H(k_0,k)$ | $\widehat{\gamma}_n^{(1.90)}(k_0,k)$ | $\widetilde{\gamma}_n^H(k_0,k)$ | $\widehat{\gamma}_n^{(1.61)}(k_0,k)$ |
|---|---|---|---|---|---|---|
| $n$ | Burr (2,1) | Burr (2,1) | Fréchet (1) | Fréchet (1) | Pareto (2) | Pareto (2) |
| 200 | 0.7020 | 0.9860 | 0.6860 | 0.9920 | 0.9370 | 0.9700 |
| 300 | 0.6930 | 0.9880 | 0.7000 | 0.9810 | 0.9310 | 0.9680 |
| 400 | 0.7030 | 0.9820 | 0.6870 | 0.9880 | 0.9140 | 0.9650 |
| 500 | 0.7165 | 0.9730 | 0.7060 | 0.9800 | 0.9250 | 0.9620 |
| 600 | 0.7085 | 0.9780 | 0.7320 | 0.9800 | 0.9060 | 0.9430 |
| 700 | 0.7100 | 0.9770 | 0.7350 | 0.9690 | 0.8980 | 0.9490 |
| 800 | 0.7130 | 0.9820 | 0.7440 | 0.9750 | 0.9260 | 0.9570 |
| 900 | 0.7115 | 0.9740 | 0.6970 | 0.9610 | 0.9020 | 0.9380 |
| 1000 | 0.7220 | 0.9745 | 0.7320 | 0.9770 | 0.9120 | 0.9400 |
| 1100 | 0.7245 | 0.9760 | 0.7170 | 0.9580 | 0.9050 | 0.9370 |
| 1200 | 0.7295 | 0.9700 | 0.7050 | 0.9640 | 0.8910 | 0.9150 |
| 1300 | 0.7255 | 0.9665 | 0.7200 | 0.9620 | 0.8940 | 0.9280 |
| 1400 | 0.7230 | 0.9615 | 0.7400 | 0.9730 | 0.9020 | 0.9320 |
| 1500 | 0.7400 | 0.9640 | 0.7070 | 0.9620 | 0.8970 | 0.9310 |
| 1600 | 0.7180 | 0.9610 | 0.7290 | 0.9700 | 0.8730 | 0.9150 |
| 1700 | 0.7455 | 0.9665 | 0.7270 | 0.9630 | 0.8760 | 0.9280 |
| 1800 | 0.7310 | 0.9615 | 0.7330 | 0.9760 | 0.8930 | 0.9110 |
| 1900 | 0.7160 | 0.9605 | 0.7350 | 0.9610 | 0.8940 | 0.9110 |
| 2000 | 0.7325 | 0.9620 | 0.7250 | 0.9600 | 0.8900 | 0.9080 |

Lastly we compare the coverage probability and the confidence length based on the asymptotic normality of the two location invariant Hill estimators. Note that Corollary 2.2 of Fraga Alves (2001) provides us the $1-\theta$ asymptotic confidence interval of $\gamma$, i.e.

$$I_N(1-\theta) = \left\{ \frac{\widetilde{\gamma}_n^H(k_0,k)}{1+z_{\theta/2}/\sqrt{k_0}}, \frac{\widetilde{\gamma}_n^H(k_0,k)}{1-z_{\theta/2}/\sqrt{k_0}} \right\}. \quad (3.2)$$

We can also construct a new asymptotic confidence interval of $\gamma$ based on Corollary 2.1 by choosing $\alpha_0$ such that $b_{\alpha_0}(\gamma) = 0$ for given $\gamma$. The asymptotic confidence interval of $\gamma$ is defined by

$$I_N(1-\theta) = \left\{ \frac{\widehat{\gamma}_n^{(\alpha_0)}(k_0,k)}{1+(V_{\alpha_0}/k_0)^{1/2} z_{\theta/2}}, \frac{\widehat{\gamma}_n^{(\alpha_0)}(k_0,k)}{1-(V_{\alpha_0}/k_0)^{1/2} z_{\theta/2}} \right\}, \quad (3.3)$$

where $k_0$ in (3.3) is the same as the one in (3.2) and $z_{\theta/2}$ is the critical value of the standard normal distribution at level $\theta/2$; that is, $1 - \Phi(z_{\theta/2}) = \theta/2$. We drew 2000 random samples from the Fréchet (1), Burr (2,1) and Pareto (2) distributions. The simulation was repeated 2000 times, we computed the coverage probabilities of $I_N(0.95)$ and the interval lengths for $n = 200, 300, \cdots, 2000$. These coverage probabilities and interval lengths are reported in Tables 2 and 3.

We may conclude that the two estimators have comparable coverage probabilities only for Pareto distribution. Generally, the proposed estimator $\widehat{\gamma}_n^{(\alpha_0)}(k_0,k)$ has better coverage probability than that of $\tilde{\gamma}_n^H(k_0,k)$. But it has a wider confidence interval length than that of $\tilde{\gamma}_n^H(k_0,k)$ especially when the sample size $n$ is small. Meanwhile, Figure 2 tells us that there exists $\alpha_0$ such that $\widehat{\gamma}_n^{(\alpha_0)}(k_0,k)$ has asymptotic null bias and better asymptotic efficiency than $\tilde{\gamma}_n^H(k_0,k)$. This new location invariant Hill type estimator may be useful in empirical analysis.

Table 3
Average length comparison for 95% confidence levels based on (3.2) and (3.3).

| $n$ | $\tilde{\gamma}_n^H(k_0,k)$ | $\widehat{\gamma}_n^{(2.37)}(k_0,k)$ | $\tilde{\gamma}_n^H(k_0,k)$ | $\widehat{\gamma}_n^{(1.90)}(k_0,k)$ | $\tilde{\gamma}_n^H(k_0,k)$ | $\widehat{\gamma}_n^{(1.61)}(k_0,k)$ |
|---|---|---|---|---|---|---|
|  | Burr (2,1) | Burr (2,1) | Fréchet (1) | Fréchet (1) | Pareto (2) | Pareto (2) |
| 200 | 0.6680 | 6.4595 | 0.6672 | 6.4281 | 1.1188 | 3.4414 |
| 300 | 0.5360 | 2.9143 | 0.5378 | 2.9436 | 0.9174 | 2.5042 |
| 400 | 0.4670 | 2.1501 | 0.4681 | 2.1292 | 0.8098 | 2.1103 |
| 500 | 0.4206 | 1.7727 | 0.4226 | 1.7903 | 0.7366 | 1.8655 |
| 600 | 0.3895 | 1.5558 | 0.3879 | 1.5590 | 0.6826 | 1.7042 |
| 700 | 0.3643 | 1.4038 | 0.3637 | 1.4130 | 0.6461 | 1.5912 |
| 800 | 0.3442 | 1.2775 | 0.3423 | 1.2887 | 0.6120 | 1.4909 |
| 900 | 0.3291 | 1.2146 | 0.3305 | 1.2292 | 0.5856 | 1.4322 |
| 1000 | 0.3148 | 1.1447 | 0.3138 | 1.1255 | 0.5616 | 1.3519 |
| 1100 | 0.3034 | 1.0865 | 0.3030 | 1.1009 | 0.5484 | 1.3237 |
| 1200 | 0.2928 | 1.0485 | 0.2948 | 1.0565 | 0.5282 | 1.2601 |
| 1300 | 0.2840 | 1.0014 | 0.2839 | 0.9973 | 0.5134 | 1.2306 |
| 1400 | 0.2754 | 0.9734 | 0.2762 | 0.9698 | 0.5013 | 1.1996 |
| 1500 | 0.2685 | 0.9418 | 0.2693 | 0.9373 | 0.4871 | 1.1606 |
| 1600 | 0.2625 | 0.9182 | 0.2618 | 0.9084 | 0.4755 | 1.1283 |
| 1700 | 0.2551 | 0.8808 | 0.2568 | 0.8926 | 0.4648 | 1.1024 |
| 1800 | 0.2505 | 0.8671 | 0.2502 | 0.8478 | 0.4572 | 1.0817 |
| 1900 | 0.2460 | 0.8453 | 0.2452 | 0.8363 | 0.4482 | 1.0618 |
| 2000 | 0.2407 | 0.8237 | 0.2409 | 0.8234 | 0.4395 | 1.0421 |

## 4. Proofs

Before proving the main results, we need some lemmas.

**Lemma 4.1.** *Suppose that (1.4) holds for $\gamma > 0$ and $\rho < 0$, then for any $\varepsilon$, $\delta > 0$, there exists $t_0 = t_0(\varepsilon, \delta)$ such that for all $t > t_0$ and $x > 1$,*

$$\left| \frac{\ln \frac{U(tx) - U(t)}{a(t)} - \ln D_\gamma(x)}{A(t)} - B_{\gamma,\rho}(x) \right| \leq \varepsilon \left(1 + \gamma \frac{x^{\gamma+\rho+\delta}}{x^\gamma - 1}\right),$$

*where $D_\gamma(x) = (x^\gamma - 1)/\gamma$ and*

$$B_{\gamma,\rho}(x) = \begin{cases} \dfrac{\gamma}{\gamma+\rho} \dfrac{x^{\gamma+\rho} - 1}{x^\gamma - 1}, & \text{if } \gamma + \rho \neq 0, \\ \dfrac{\gamma}{x^\gamma - 1} \ln x, & \text{if } \gamma + \rho = 0. \end{cases}$$

*Proof.* We only consider the case of $\gamma + \rho \neq 0$. Note that by (1.4) we have

$$\frac{\frac{U(tx)-U(t)}{a(t)D_\gamma(x)} - 1}{A(t)} \to \frac{\gamma}{\gamma + \rho} \frac{x^{\gamma+\rho} - 1}{x^\gamma - 1}.$$

This implies

$$\lim_{t \to \infty} \frac{\ln \frac{U(tx)-U(t)}{a(t)} - \ln D_\gamma(x)}{A(t)} = B_{\gamma,\rho}(x) \tag{4.1}$$

since $\ln(1+x) \sim x$ as $x \to 0$. By Theorem 2.3.6 in de Haan and Ferreira (2006), for any $\varepsilon, \delta > 0$, there exists $t_0 = t_0(\varepsilon, \delta)$ such that for all $t > t_0$ and $x > 1$,

$$\left| \frac{\frac{U(tx)-U(t)}{a(t)} - D_\gamma(x)}{A(t)} - \Psi_{\gamma,\rho}(x) \right| \leq \varepsilon x^{\gamma+\rho+\delta}.$$

Hence,

$$\left| \frac{\ln \frac{U(tx)-U(t)}{a(t)} - \ln D_\gamma(x)}{A(t)} - B_{\gamma,\rho}(x) \right|$$

$$= \left| \frac{\ln \frac{U(tx)-U(t)}{a(t)D_\gamma(x)} - \left( \frac{U(tx)-U(t)}{a(t)D_\gamma(x)} - 1 \right)}{A(t)} \right|$$

$$+ \frac{1}{D_\gamma(x)} \left| \frac{\frac{U(tx)-U(t)}{a(t)} - D_\gamma(x)}{A(t)} - \Psi_{\gamma,\rho}(x) \right|$$

$$\leq \varepsilon \left( 1 + \gamma \frac{x^{\gamma+\rho+\delta}}{x^\gamma - 1} \right).$$

The proof is complete. □

**Lemma 4.2.** *Under the condition of Lemma 4.1, for $x > 1$, $y > 1$ we have*

$$\lim_{t \to \infty} \frac{\ln \frac{U(tx)-U(t)}{U(ty)-U(t)} - \ln \frac{x^\gamma-1}{y^\gamma-1}}{A(t)} = F_{\gamma,\rho}(x, y).$$

*Moreover, for any $\varepsilon, \delta > 0$, there exists $t_0 = t_0(\varepsilon, \delta)$ such that for all $t > t_0$ and $x > y > 1$,*

$$\left| \frac{\ln \frac{U(tx)-U(t)}{U(ty)-U(t)} - \ln \frac{D_\gamma(x)}{D_\gamma(y)}}{A(t)} - F_{\gamma,\rho}(x, y) \right| \leq T_{\gamma,\rho}(x, y),$$

*where*

$$F_{\gamma,\rho}(x, y) = \begin{cases} \frac{\gamma}{\gamma + \rho} \left( \frac{x^{\gamma+\rho} - 1}{x^\gamma - 1} - \frac{y^{\gamma+\rho} - 1}{y^\gamma - 1} \right), & \text{if } \gamma + \rho \neq 0, \\ \frac{\gamma}{x^\gamma - 1} \ln x - \frac{\gamma}{y^\gamma - 1} \ln y, & \text{if } \gamma + \rho = 0 \end{cases}$$

*and $T_{\gamma,\rho}(x, y) = \varepsilon(2 + \gamma \{x^{\gamma+\rho+\delta} + y^{\gamma+\rho+\delta}\}/\{y^\gamma - 1\})$.*

*Proof.* We only consider the case of $\gamma + \rho \neq 0$, similar arguments for $\gamma + \rho = 0$. Note that (4.1) implies

$$\lim_{t\to\infty} \frac{\ln\frac{U(tx)-U(t)}{U(ty)-U(t)} - \ln\frac{x^\gamma-1}{y^\gamma-1}}{A(t)} = \frac{\gamma}{\gamma+\rho}\left(\frac{x^{\gamma+\rho}-1}{x^\gamma-1} - \frac{y^{\gamma+\rho}-1}{y^\gamma-1}\right)$$

for $x > 1$, $y > 1$. For any $\varepsilon, \delta > 0$, there exists $t_0 = t_0(\varepsilon,\delta)$ such that for all $t > t_0$ and $x > y > 1$,

$$\left|\frac{\ln\frac{U(tx)-U(t)}{U(ty)-U(t)} - \ln\frac{D_\gamma(x)}{D_\gamma(y)}}{A(t)} - F_{\gamma,\rho}(x,y)\right|$$

$$= \left|\frac{\ln\frac{U(tx)-U(t)}{a(t)} - \ln D_\gamma(x)}{A(t)} - \frac{\ln\frac{U(ty)-U(t)}{a(t)} - \ln D_\gamma(y)}{A(t)} - F_{\gamma,\rho}(x,y)\right|$$

$$\leq \left|\frac{\ln\frac{U(tx)-U(t)}{a(t)} - \ln D_\gamma(x)}{A(t)} - B_{\gamma,\rho}(x)\right|$$

$$+ \left|\frac{\ln\frac{U(ty)-U(t)}{a(t)} - \ln D_\gamma(y)}{A(t)} - B_{\gamma,\rho}(y)\right|$$

$$+ \left|B_{\gamma,\rho}(x) - B_{\gamma,\rho}(y) - F_{\gamma,\rho}(x,y)\right|.$$

By using Lemma 4.1, we get the desired result. □

*Proof of Theorem 2.1.* For the asymptotic distributional representation of $M_n^{(\alpha)}(k_0,k)$, we only consider the case of $\gamma + \rho \neq 0$. By Lemma 4.2, we have

$$\ln\frac{U(tx)-U(t)}{U(ty)-U(t)}$$

$$= \ln\frac{x^\gamma-1}{y^\gamma-1} + \frac{\gamma}{\gamma+\rho}\left[\frac{x^{\gamma+\rho}-1}{x^\gamma-1} - \frac{y^{\gamma+\rho}-1}{y^\gamma-1}\right]A(t)(1+o(1))$$

$$= \gamma\ln\frac{x}{y} + \ln\left(1 + \frac{y^{-\gamma}-x^{-\gamma}}{1-y^{-\gamma}}\right)$$

$$+ \frac{\gamma}{\gamma+\rho}\left[\frac{x^\rho-x^{-\gamma}}{1-x^{-\gamma}} - \frac{y^\rho-y^{-\gamma}}{1-y^{-\gamma}}\right]A(t)(1+o(1)).$$

Let $Y_1, Y_2, \cdots, Y_n$ be i.i.d Pareto r.v.s with $F_Y(y) = 1 - 1/y$, $y \geq 1$ and let $Y_{1,n} \leq Y_{2,n} \leq \cdots \leq Y_{n,n}$ denote the order statistics of $Y_1, Y_2, \cdots, Y_n$. Now replace $t$ by $Y_{n-k,n}$, $x$ by $Y_{n-i,n}/Y_{n-k,n}$ and $y$ by $Y_{n-k_0,n}/Y_{n-k,n}$, respectively, and note that $\{X_i\}_{i=1}^n \stackrel{d}{=} \{U(Y_i)\}_{i=1}^n$, $(Y_{n-i,n}/Y_{n-k,n})^{-\gamma} < (Y_{n-k_0,n}/Y_{n-k,n})^{-\gamma} \to 0$ in probability uniformly for $i = 0, 1, \cdots, k_0 - 1$, and $\{Y_{n-i,n}/Y_{n-k_0,n}\}_{i=1}^{k_0-1} \stackrel{d}{=} \{Y_{k_0-i,k_0}\}_{i=1}^{k_0-1}$. So,

$$\ln \frac{X_{n-i,n} - X_{n-k,n}}{X_{n-k_0,n} - X_{n-k,n}}$$

$$\stackrel{d}{=} \ln \frac{U(Y_{n-i,n}) - U(Y_{n-k,n})}{U(Y_{n-k_0,n}) - U(Y_{n-k,n})}$$

$$= \gamma \ln \frac{Y_{n-i,n}}{Y_{n-k_0,n}} + \left(\frac{Y_{n-k_0,n}}{Y_{n-k,n}}\right)^{-\gamma} \left[1 - \left(\frac{Y_{n-i,n}}{Y_{n-k_0,n}}\right)^{-\gamma}\right](1+o_P(1))$$

$$+ \frac{\gamma}{\gamma + \rho}\left[\left(\frac{Y_{n-k_0,n}}{Y_{n-k,n}}\right)^{-\gamma}\left(1 - \left(\frac{Y_{n-i,n}}{Y_{n-k_0,n}}\right)^{-\gamma}\right)\right.$$

$$\left. - \left(\frac{Y_{n-k_0,n}}{Y_{n-k,n}}\right)^{\rho}\left(1 - \left(\frac{Y_{n-i,n}}{Y_{n-k_0,n}}\right)^{\rho}\right)\right] A(Y_{n-k,n})(1+o_P(1))$$

$$\stackrel{d}{=} \gamma \ln Y_{k_0-i,k_0} + \left[1 - Y_{k_0-i,k_0}^{-\gamma}\right]\left(\frac{k_0}{k}\right)^{\gamma}(1+o_P(1))$$

$$+ \frac{\gamma}{\gamma + \rho}\left[\left(1 - Y_{k_0-i,k_0}^{-\gamma}\right)\left(\frac{k_0}{k}\right)^{\gamma}\right.$$

$$\left. - \left(1 - Y_{k_0-i,k_0}^{\rho}\right)\left(\frac{k_0}{k}\right)^{-\rho}\right] A(n/k)(1+o_P(1)).$$

By Taylor's expansion, we may get

$$M_n^{(\alpha)}(k_0, k)$$

$$\stackrel{d}{=} \frac{1}{k_0} \sum_{i=1}^{k_0-1} \left\{ (\gamma \ln Y_{k_0-i,k_0})^{\alpha} \right.$$

$$+ \alpha \gamma^{\alpha} (\ln Y_{k_0-i,k_0})^{\alpha-1} \frac{1 - Y_{k_0-i,k_0}^{-\gamma}}{\gamma}\left(\frac{k_0}{k}\right)^{\gamma}(1+o_P(1))$$

$$+ \frac{\alpha \gamma^{\alpha}}{\gamma + \rho}(\ln Y_{k_0-i,k_0})^{\alpha-1}\left[\left(1 - Y_{k_0-i,k_0}^{-\gamma}\right)\left(\frac{k_0}{k}\right)^{\gamma}\right.$$

$$\left.\left. - (1 - Y_{k_0-i,k_0}^{\rho})\left(\frac{k_0}{k}\right)^{-\rho}\right] A(n/k)(1+o_P(1)) \right\}$$

$$+ A^2(n/k) O_P\left(\left(\frac{k_0}{k}\right)^{2\gamma} + \left(\frac{k_0}{k}\right)^{-2\rho} + \left(\frac{k_0}{k}\right)^{\gamma-\rho}\right).$$

Note that $E[(\ln Y_1)^{\alpha}] = \Gamma(\alpha + 1)$, $Var[(\ln Y_1)^{\alpha}] = \sigma_{\alpha}$ and

$$E\left[(\ln Y_1)^{\alpha-1}(Y_1^{\rho} - 1)/\rho\right] = \mu_{\alpha}(\rho).$$

So,

$$P_n^{(\alpha)} = \frac{\frac{1}{k_0}\sum_{i=1}^{k_0}(\ln Y_i)^{\alpha} - \Gamma(\alpha+1)}{\sigma_{\alpha}/\sqrt{k_0}} \stackrel{d}{\to} N(0,1). \tag{4.2}$$

By (4.2) and the law of large numbers, the asymptotic distributional representation of $M_n^{(\alpha)}(k_0, k)$ is

$$M_n^{(\alpha)}(k_0, k) \stackrel{d}{=} \Gamma(\alpha+1)\gamma^\alpha + \gamma^\alpha \frac{\sigma_\alpha}{\sqrt{k_0}} P_n^{(\alpha)} + \alpha\gamma^\alpha \mu_\alpha(-\gamma) \left(\frac{k_0}{k}\right)^\gamma (1+o_P(1))$$
$$+ \frac{\alpha\gamma^\alpha \rho}{\gamma+\rho} \mu_\alpha(\rho) A(n/k) \left(\frac{k_0}{k}\right)^{-\rho} (1+o_P(1)), \quad (4.3)$$

which is the desired result. □

*Proof of Theorem 2.2.* Note that both $\{X_i\}_{i=1}^n \stackrel{d}{=} \{U(Y_i)\}_{i=1}^n$ and $\{Y_{n-i,n} / Y_{n-k_0,n}\}_{i=1}^{k_0-1} \stackrel{d}{=} \{Y_{k_0-i,k_0}\}_{i=1}^{k_0-1}$ hold. For $\gamma + \rho \neq 0$, by using the conditions of Theorem 2.2, following joint expansion in distribution,

$$\left(\left(\frac{M_n^{(2\alpha)}(k_0, k)}{\Gamma(2\alpha+1)}\right)^{1/2}, \frac{\Gamma(\alpha)}{M_n^{(\alpha-1)}(k_0, k)}\right) \stackrel{d}{=} (W_n, Z_n) \quad (4.4)$$

holds, where

$$W_n = \gamma^\alpha \left[1 + \left(\frac{\sigma_{2\alpha}}{2\Gamma(2\alpha+1)} \frac{P_n^{(2\alpha)}}{\sqrt{k_0}} + \frac{\mu_{2\alpha}(-\gamma)}{2\Gamma(2\alpha)} \left(\frac{k_0}{k}\right)^\gamma \right. \right.$$
$$\left.\left. + \frac{\rho\mu_{2\alpha}(\rho)}{2(\gamma+\rho)\Gamma(2\alpha)} A\left(\frac{n}{k}\right) \left(\frac{k_0}{k}\right)^{-\rho}\right)(1+o_P(1))\right]$$

and

$$Z_n = \gamma^{1-\alpha} \left[1 - \left(\frac{\sigma_{\alpha-1}}{\Gamma(\alpha)} \frac{P_n^{(\alpha-1)}}{\sqrt{k_0}} + \frac{\mu_{\alpha-1}(-\gamma)}{\Gamma(\alpha-1)} \left(\frac{k_0}{k}\right)^\gamma \right.\right.$$
$$\left.\left. + \frac{\rho\mu_{\alpha-1}(\rho)}{(\gamma+\rho)\Gamma(\alpha-1)} A\left(\frac{n}{k}\right) \left(\frac{k_0}{k}\right)^{-\rho}\right)(1+o_P(1))\right].$$

Note that (4.4) may be deduced from Wold device and Delta-method, as for arbitrary $a, b \in R$, we have

$$a\left(M_n^{(2\alpha)}(k_0,k)/\Gamma(2\alpha+1)\right)^{1/2} + b\Gamma(\alpha)\left(M_n^{(\alpha-1)}(k_0,k)\right)^{-1}$$
$$\stackrel{d}{=} a\left(\frac{1}{\Gamma(2\alpha+1)}\frac{1}{k_0}\sum_{i=1}^{k_0-1}\left(\ln\frac{U(Y_{n-i,n}) - U(Y_{n-k,n})}{U(Y_{n-k_0,n}) - U(Y_{n-k,n})}\right)^{2\alpha}\right)^{1/2}$$
$$+ b\Gamma(\alpha)\left(\frac{1}{k_0}\sum_{i=1}^{k_0-1}\left(\ln\frac{U(Y_{n-i,n}) - U(Y_{n-k,n})}{U(Y_{n-k_0,n}) - U(Y_{n-k,n})}\right)^{\alpha-1}\right)^{-1}$$
$$\stackrel{d}{=} a\left(\frac{Q_n^{(2\alpha)}}{\Gamma(2\alpha+1)}\right)^{1/2} + b\Gamma(\alpha)\left(Q_n^{(\alpha-1)}\right)^{-1}$$
$$= aW_n + b\gamma^{1-\alpha}\frac{1}{1+(1-Z_n\gamma^{\alpha-1})}$$
$$= aW_n + bZ_n.$$

The last step follows since $Z_n \xrightarrow{P} \gamma^{1-\alpha}$ and $h(1 - Z_n\gamma^{\alpha-1}) - h(0) = h'(0)(1 + o_P(1))(1 - Z_n\gamma^{\alpha-1})$ with $h(x) = 1/(1+x)$, where

$$\begin{aligned} Q_n^{(\beta)} &= \frac{1}{k_0} \sum_{i=1}^{k_0-1} \Bigg\{ (\gamma \ln Y_{k_0-i,k_0})^\beta \\ &\quad + \beta\gamma^\beta (\ln Y_{k_0-i,k_0})^{\beta-1} \frac{1 - Y_{k_0-i,k_0}^{-\gamma}}{\gamma} \left(\frac{k_0}{k}\right)^\gamma (1 + o_P(1)) \\ &\quad + \frac{\beta\gamma^\beta}{\gamma+\rho}(\ln Y_{k_0-i,k_0})^{\beta-1} \Bigg[ \left(1 - Y_{k_0-i,k_0}^{-\gamma}\right)\left(\frac{k_0}{k}\right)^\gamma \\ &\quad - \left(1 - Y_{k_0-i,k_0}^{\rho}\right)\left(\frac{k_0}{k}\right)^{-\rho} \Bigg] A(n/k)(1+o_P(1)) \Bigg\} \\ &\quad + A^2(n/k) O_P\left(\left(\frac{k_0}{k}\right)^{2\gamma} + \left(\frac{k_0}{k}\right)^{-2\rho} + \left(\frac{k_0}{k}\right)^{\gamma-\rho}\right) \end{aligned}$$

for $\beta \geq 1$. So,

$$\begin{aligned} \widehat{\gamma}_n^{(\alpha)}(k_0, k) &\stackrel{d}{=} Z_n W_n \\ &= \gamma \Bigg[ 1 + \frac{1}{\sqrt{k_0}}\left(\frac{\sigma_{2\alpha}}{2\Gamma(2\alpha+1)} P_n^{(2\alpha)} - \frac{\sigma_{\alpha-1}}{\Gamma(\alpha)} P_n^{(\alpha-1)}\right) \\ &\quad + \left(\frac{\mu_{2\alpha}(-\gamma)}{2\Gamma(2\alpha)} - \frac{\mu_{\alpha-1}(-\gamma)}{\Gamma(\alpha-1)}\right)\left(\frac{k_0}{k}\right)^\gamma \\ &\quad + \frac{\rho}{\gamma+\rho}\left(\frac{\mu_{2\alpha}(\rho)}{2\Gamma(2\alpha)} - \frac{\mu_{\alpha-1}(\rho)}{\Gamma(\alpha-1)}\right) A\left(\frac{n}{k}\right)\left(\frac{k_0}{k}\right)^{-\rho}(1+o_P(1)) \Bigg] \\ &\quad + o_P\left(\frac{1}{\sqrt{k_0}}\right) + o_P\left(\left(\frac{k_0}{k}\right)^\gamma\right). \end{aligned}$$

Denote

$$Q_n^{(\alpha)} = \frac{\sigma_{2\alpha}}{2\Gamma(2\alpha+1)} P_n^{(2\alpha)} - \frac{\sigma_{\alpha-1}}{\Gamma(\alpha)} P_n^{(\alpha-1)}, \quad T_n^{(\alpha)} = \frac{Q_n^{(\alpha)}}{\sqrt{Var(Q_n^{(\alpha)})}}.$$

Let $f_{k_0}(t)$ denote the characteristic function of $Q_n^{(\alpha)}$. Noting the expression of $P_n^{(\alpha)}$ in (4.2) and that of $V_\alpha$ in (2.2), we have

$$\begin{aligned} f_{k_0}(t) &= E\exp\left\{it\frac{\sigma_{2\alpha}}{2\Gamma(2\alpha+1)} P_n^{(2\alpha)} - it\frac{\sigma_{\alpha-1}}{\Gamma(\alpha)} P_n^{(\alpha-1)}\right\} \\ &= E\exp\left\{\frac{it}{\sqrt{k_0}}\sum_{j=1}^{k_0} \frac{(\ln Y_j)^{2\alpha} - \Gamma(2\alpha+1)}{2\Gamma(2\alpha+1)} \right. \\ &\quad \left. - \frac{it}{\sqrt{k_0}}\sum_{j=1}^{k_0} \frac{(\ln Y_j)^{\alpha-1} - \Gamma(\alpha)}{\Gamma(\alpha)}\right\} \end{aligned}$$

$$= \prod_{j=1}^{k_0} E \exp\left\{\frac{it}{\sqrt{k_0}}\left(\frac{(\ln Y_j)^{2\alpha} - \Gamma(2\alpha+1)}{2\Gamma(2\alpha+1)} - \frac{(\ln Y_j)^{\alpha-1} - \Gamma(\alpha)}{\Gamma(\alpha)}\right)\right\}$$

$$= \prod_{j=1}^{k_0}\left\{1 - \frac{t^2}{2k_0}E\left(\frac{(\ln Y_j)^{2\alpha} - \Gamma(2\alpha+1)}{2\Gamma(2\alpha+1)} - \frac{(\ln Y_j)^{\alpha-1} - \Gamma(\alpha)}{\Gamma(\alpha)}\right)^2 + o\left(\frac{1}{k_0}\right)\right\}$$

$$= \left\{1 - \frac{t^2}{2k_0}V_\alpha + o\left(\frac{1}{k_0}\right)\right\}^{k_0} \to \exp(-t^2 V_\alpha/2).$$

So, $Q_n^{(\alpha)}$ converges in distribution to a normal r.v. with null bias and variance $V_\alpha$ and $T_n^{(\alpha)}$ is an asymptotically standard normal r.v. Noting $\gamma[\mu_{2\alpha}(-\gamma)/\{2\Gamma(2\alpha)\} - \mu_{\alpha-1}(-\gamma)/\Gamma(\alpha-1)] = b_\alpha(\gamma)$ and $\{\gamma\rho/(\gamma+\rho)\}[\mu_{2\alpha}(\rho)/\{2\Gamma(2\alpha)\} - \mu_{\alpha-1}(\rho)/\Gamma(\alpha-1)] = -\{\gamma/(\gamma+\rho)\}b_\alpha(-\rho)$, we obtain the asymptotic distributional representation of $\widehat{\gamma}_n^{(\alpha)}(k_0,k)$. The remaining part of Theorem 2.2 is immediate. □

*Proof of Theorem 2.3.* Assume that $A(t) \sim ct^\rho$ with $c \neq 0$ and $\rho < 0$, then according to the range of the pair $(\gamma, \rho)$, we need to consider the following two cases.

(i). <u>The $\gamma \leq -\rho$ case</u>. Combining (2.4) and Theorem 2.2, we have $R_n = o((k_0/k)^\gamma)$ and

$$\widehat{\gamma}_n^{(\alpha)}(k_0,k) \stackrel{d}{=} \gamma + \frac{\gamma\sqrt{V_\alpha}}{\sqrt{k_0}}T_n^{(\alpha)} + b_\alpha(\gamma)\left(\frac{k_0}{k}\right)^\gamma + o_P\left(\frac{1}{\sqrt{k_0}}\right) + o_P\left(\left(\frac{k_0}{k}\right)^\gamma\right)$$

and

$$MSE_\infty\left(\widehat{\gamma}_n^{(\alpha)}(k_0,k)\right) = \frac{\gamma^2 V_\alpha}{k_0} + b_\alpha^2(\gamma)\left(\frac{k_0}{k}\right)^{2\gamma}.$$

So, the sequence of $k_0(n)$ that minimizes $MSE_\infty(\widehat{\gamma}_n^{(\alpha)}(k_0,k))$ is

$$k_0^{opt} = k_0^{(1)} = \left[\frac{\gamma V_\alpha}{2b_\alpha^2(\gamma)}\right]^{\frac{1}{2\gamma+1}} \cdot k^{\frac{2\gamma}{2\gamma+1}}.$$

(ii). <u>The $\gamma > -\rho$ case</u>. Noting the expression of $R_n$ in (2.4), we have $R_n - c\gamma b_\alpha(-\rho)/(\gamma+\rho)(k_0/n)^{-\rho}(1+o_P(1))$ and

$$\widehat{\gamma}_n^{(\alpha)}(k_0,k) \stackrel{d}{=} \gamma + \frac{\gamma\sqrt{V_\alpha}}{\sqrt{k_0}}T_n^{(\alpha)} + b_\alpha(\gamma)\left(\frac{k_0}{k}\right)^\gamma - \frac{c\gamma b_\alpha(-\rho)}{\gamma+\rho}\left(\frac{k_0}{n}\right)^{-\rho}$$
$$+ o_P\left(\frac{1}{\sqrt{k_0}}\right) + o_P\left(\left(\frac{k_0}{k}\right)^\gamma\right) + o_P\left(\left(\frac{k_0}{n}\right)^{-\rho}\right). \quad (4.5)$$

We need to investigate the related weights in (4.5). Consider any sequence $k_0$ satisfying $k_0 = O(k^{\gamma/(\gamma+\rho)} \cdot n^{\rho/(\gamma+\rho)})$, then $(k_0/k)^\gamma = O((k_0/n)^{-\rho})$. Moreover, if $k_0 \ll k^{\gamma/(\gamma+\rho)} \cdot n^{\rho/(\gamma+\rho)}$, then $(k_0/k)^\gamma \ll (k_0/n)^{-\rho}$; if $k_0 \gg k^{\gamma/(\gamma+\rho)} \cdot n^{\rho/(\gamma+\rho)}$, then $(k_0/k)^\gamma \gg (k_0/n)^{-\rho}$. In order to obtain the optimal choice of $k_0$, we have to consider the relationship of $k$ and $n$.

(a). Firstly, suppose $k^{2\gamma/(2\gamma+1)} \ll n^{-2\rho/(-2\rho+1)}$, then

$$k^{\frac{\gamma}{\gamma+\rho}} \cdot n^{\frac{\rho}{\gamma+\rho}} \ll k^{\frac{2\gamma}{2\gamma+1}} \left(\ll n^{\frac{-2\rho}{-2\rho+1}}\right).$$

If the sequence $k_0$ satisfies $k_0 = O(k^{2\gamma/(2\gamma+1)})$, then $1/\sqrt{k_0} = O((k_0/k)^\gamma)$. Moreover, if $k_0 \ll k^{2\gamma/(2\gamma+1)}$, then $1/\sqrt{k_0} \gg (k_0/k)^\gamma$, $1/\sqrt{k_0} \gg (k_0/k)^{-\rho}$, and $MSE_\infty$ $(\widehat{\gamma}_n^{(\alpha)}(k_0, k)) = \gamma^2 V_\alpha/k_0$ is a decreasing function of $k_0$; If $k_0 \gg k^{2\gamma/(2\gamma+1)}$, then $(k_0/k)^\gamma \gg 1/\sqrt{k_0}$, $(k_0/k)^\gamma \gg (k_0/n)^{-\rho}$ and $MSE_\infty$ $(\widehat{\gamma}_n^{(\alpha)}(k_0, k)) = b_\alpha^2(\gamma)$ $(k_0/k)^{2\gamma}$ is an increasing function of $k_0$. So, we choose $k_0 = O(k^{2\gamma/(2\gamma+1)})$ in order to balance the bias and variance of the estimator $\widehat{\gamma}_n^{(\alpha)}(k_0, k)$, hence $MSE_\infty(\widehat{\gamma}_n^{(\alpha)}(k_0, k)) = \gamma^2 V_\alpha/k_0 + b_\alpha^2(\gamma)(k_0/k)^{2\gamma}$ and

$$k_0^{opt} \equiv k_0^{(1)} = \left[\frac{\gamma V_\alpha}{2b_\alpha^2(\gamma)}\right]^{\frac{1}{2\gamma+1}} \cdot k^{\frac{2\gamma}{2\gamma+1}}.$$

(b). If $k^{2\gamma/(2\gamma+1)} \gg n^{-2\rho/(-2\rho+1)}$, we have

$$k^{\frac{\gamma}{\gamma+\rho}} \cdot n^{\frac{\rho}{\gamma+\rho}} \gg k^{\frac{2\gamma}{2\gamma+1}} \left(\gg n^{\frac{-2\rho}{-2\rho+1}}\right).$$

In this situation, if $k_0 = O(n^{-2\rho/(-2\rho+1)})$, then $1/\sqrt{k_0} = O((k_0/n)^{-\rho})$. Moreover, if $k_0 \ll n^{-2\rho/(-2\rho+1)}$, then $1/\sqrt{k_0} \gg (k_0/n)^{-\rho} \gg (k_0/k)^\gamma$ and $MSE_\infty$ $(\widehat{\gamma}_n^{(\alpha)}(k_0, k)) = \gamma^2 V_\alpha/k_0$; If $n^{-2\rho/(-2\rho+1)} \ll k_0 \ll k^{\gamma/(\gamma+\rho)} \cdot n^{\rho/(\gamma+\rho)}$, then $(k_0/n)^{-\rho} \gg 1/\sqrt{k_0}$, $(k_0/n)^{-\rho} \gg (k_0/k)^\gamma$ and $MSE_\infty$ $(\widehat{\gamma}_n^{(\alpha)}(k_0, k)) = [c\gamma\, b_\alpha(-\rho)/(\gamma+\rho)]^2$ $(k_0/n)^{-2\rho}$; if $k_0 \gg k^{\gamma/(\gamma+\rho)} \cdot n^{\rho/(\gamma+\rho)}$, then $(k_0/k)^\gamma \gg (k_0/n)^{-\rho} \gg 1/\sqrt{k_0}$.

Choosing the sequence $k_0 = O(n^{-2\rho/(-2\rho+1)})$ to balance the bias and variance is possible. Then

$$MSE_\infty\left(\widehat{\gamma}_n^{(\alpha)}(k_0, k)\right) = \frac{\gamma^2 V_\alpha}{k_0} + \left[\frac{c\gamma b_\alpha(-\rho)}{\gamma+\rho}\right]^2 \left(\frac{k_0}{n}\right)^{-2\rho}.$$

We get the locally minimized MSE for the sequence

$$k_0^{opt} = k_0^{(2)} = \left[\frac{(\gamma+\rho)^2 V_\alpha}{-2\rho c^2 b_\alpha^2(-\rho)}\right]^{\frac{1}{-2\rho+1}} \cdot n^{\frac{-2\rho}{-2\rho+1}}.$$

But if any sequence $k_0 = O(k^{\gamma/(\gamma+\rho)} \cdot n^{\rho/(\gamma+\rho)})$ and

$$MSE_\infty(\widehat{\gamma}_n^{(\alpha)}(k_0, k)) = \left[b_\alpha(\gamma)\left(\frac{k_0}{k}\right)^\gamma - \frac{c\gamma b_\alpha(-\rho)}{\gamma+\rho}\left(\frac{k_0}{n}\right)^{-\rho}\right]^2, \qquad (4.6)$$

we need to see the sign of $cb_\alpha(\gamma)b_\alpha(-\rho)$. Note that if $cb_\alpha(\gamma)b_\alpha(-\rho) < 0$, (4.6) is an increasing function of $k_0$, and $k_0^{opt} = k_0^{(2)}$ is the solution to the optimization problem. On the other hand, suppose $cb_\alpha(\gamma)b_\alpha(-\rho) > 0$, let $k_0^{opt} = k_0^{(3)}$ be such that $MSE_\infty(\widehat{\gamma}_n^{(\alpha)}(k_0, k)) = 0$ in (4.6). Then, the optimal sample fraction is

$$k_0^{(3)} = \left[\frac{c\gamma b_\alpha(-\rho)}{(\gamma + \rho)b_\alpha(\gamma)}\right]^{\frac{1}{\gamma+\rho}} \cdot k^{\frac{\gamma}{\gamma+\rho}} \cdot n^{\frac{\rho}{\gamma+\rho}}.$$

(c). Lastly, we consider $k \sim Dn^{-\rho(2\gamma+1)/\{\gamma(-2\rho+1)\}}$ with $D \neq 0$, i.e. $k^{2\gamma/(2\gamma+1)} = O(n^{-2\rho/(-2\rho+1)})$. We derive that $k^{\gamma/(\gamma+\rho)} \cdot n^{\rho/(\gamma+\rho)}$ is of the same order of either $k^{2\gamma/(2\gamma+1)}$ or $n^{-2\rho/(-2\rho+1)}$, and $1/\sqrt{k_0}$ is of the same order of either $(k_0/k)^\gamma$ or $(k_0/n)^{-\rho}$. Hence, $k_0^{opt}$ must be the sequence such that

$$MSE_\infty(\widehat{\gamma}_n^{(\alpha)}(k_0, k)) = \frac{\gamma^2 V_\alpha}{k_0} + \left[b_\alpha(\gamma)\left(\frac{k_0}{k}\right)^\gamma - \frac{c\gamma b_\alpha(-\rho)}{\gamma+\rho}\left(\frac{k_0}{n}\right)^{-\rho}\right]^2$$

attains its minimum, which enables us to identify $k_0^{opt} \sim D_1 n^{-2\rho/(-2\rho+1)}$ with $D_1 \equiv D_1(\gamma, \rho, \alpha, c, D)$ satisfying

$$a_1 D_1^{2\gamma+1} + a_2 D_1^{\gamma-\rho+1} + a_3 D_1^{-2\rho+1} = \gamma^2 V_\alpha,$$

where $a_1 = 2\gamma b_\alpha^2(\gamma)D^{-2\gamma}$, $a_2 = [2c\gamma(\rho-\gamma)/(\gamma+\rho)]b_\alpha(\gamma)b_\alpha(-\rho)D^{-\gamma}$ and $a_3 = -2\rho[c\gamma b_\alpha(-\rho)/(\gamma+\rho)]^2$. □

## Acknowledgments

The authors would like to thank the Editor and the referee for carefully reading the paper and for their comments which greatly improved the paper.